\documentclass{article}
\usepackage[all,2cell]{xy}
\usepackage{latexsym}

\author{Fran\c cois M\'etayer\thanks{\'Equipe PPS, 
Universit\'e Paris 7-CNRS, 2 pl.\ Jussieu,
Case 7014
F75251 Paris Cedex 05 Email:
{\tt metayer@logique.jussieu.fr}}}

\newtheorem{thm}{Theorem}

\newtheorem{prop}{Proposition}

\newtheorem{lem}{Lemma}

\newcommand{\pre}{{\bf Proof.}\ }
\newcommand{\erp}{\hfill $\Box$ \par\bigskip}

\newcommand{\ctg}[1]{\hbox{\bf #1}}

\newcommand{\cb}{\circ}

\newcommand{\up}[2]{#1^{#2}}
\newcommand{\dn}[2]{#1_{#2}}
\newcommand{\x}[3]{#1_{#2}^{#3}}

\newcommand{\doubl}[2]{\ar@<2pt>[l]^{#2}\ar@<-2pt>[l]_{#1}}
\newcommand{\doubr}[2]{\ar@<2pt>[r]^{#1}\ar@<-2pt>[r]_{#2}}
\newcommand{\doubld}[2]{\ar@<2pt>[ld]^{#2}\ar@<-2pt>[ld]_{#1}}
\newcommand{\doubd}[2]{\ar@<2pt>[d]^{#2}\ar@<-2pt>[d]_{#1}}

\newcommand{\abst}[2]{\lambda #1#2}

\newcommand{\eval}{\epsilon}

\newcommand{\term}[1]{\hbox{\tt #1}}

\newcommand{\id}{\mathop {\mathrm{id}}\nolimits}

\newcommand{\Hom}[3]{{\mathrm{Hom}}_#1(#2,#3)}

\newcommand{\abs}[1]{\left | #1 \right |}

\newcommand{\ens}[1]{\{#1\}}

\newcommand{\pair}[2]{\left\langle #1,#2\right\rangle}

\title{State monads and their algebras}
\date{}

\begin{document}

\maketitle

\begin{abstract}
State monads in cartesian closed categories are those defined by
the familiar adjunction between product and exponential. We
investigate the structure of their algebras, and show that
the exponential functor is monadic provided the base category
is sufficiently regular, and the exponent is a non-empty object.
\end{abstract}

\section{Introduction}\label{sec:intro}

Let $\ctg{C}$ be a cartesian closed category. Any object $S$  in $\ctg{C}$
gives rise to a monad $T=UF$, where
\begin{displaymath}
  U: X\mapsto X^S
\end{displaymath}
and
\begin{displaymath}
  F: X\mapsto S\times X
\end{displaymath}
$T$, among other monads, has been extensively studied by 
researchers in functional programming, in order to model
computational effects (see for example 
\cite{moggi:comlcm,moggi:notcam,wadler:common}),
and called
a {\em state monad} since then. A main source of the present
work is \cite{plotkinpower:notcdm}, which investigates
the algebras of state monads, in a significantly different
setting however, and provides an explicit
equational presentation of these algebras. 

This raises the question as whether the original
functor $U$ itself is monadic: it turns out that the answer
is positive, provided $\ctg{C}$ is sufficiently regular, and
$S$ is non-empty (Theorem~\ref{thm:monadicity}). 

\subsection{Notations}
Throughout this article, $S$ will denote
the fixed object of $\ctg{C}$, on which $T$ is built. The natural
isomorphism
\begin{displaymath}
  \Hom{C}{S\times X}{Y} \to \Hom{C}{X}{Y^S}
\end{displaymath}
will be denoted by
\begin{displaymath}
  f\mapsto \up f*
\end{displaymath}
and its inverse by
\begin{displaymath}
  g\mapsto \dn g*
\end{displaymath}
In particular, for each $X$, the identity
\begin{displaymath}
  \id_{X^S}:\up XS\to\up XS
\end{displaymath}
gives rise to
\begin{displaymath}
  \dn{({\id}_{X^S})}*:S\times\up XS\to X
\end{displaymath}
which is the {\em evaluation morphism}
$\pair{s}{f}\mapsto f[s]$, and will be denoted by $\eval_X$.
It is of course the counit of the adjunction.
Finally, for each object $X$ in $\ctg{C}$, the projections
$S\times X\to S$ and $S\times X\to X$ will be denoted by
$\dn pX$, $\dn qX$ respectively. Note that $p$,
$q$ are natural transformations $F\to 1$.

\section{The case of $\ctg{Sets}$}\label{sec:sets}

We suppose in this section that $\ctg{C}$ is the category of sets.
Let us first recall the concrete meaning of $(T,\mu,\eta)$ in this case.
For each set $X$, $TX=(S\times X)^S$, so that the unit
$\dn{\eta}{X}$ associates to each $x\in X$ the map
$\abst{s}{\pair{s}{x}}$, an element of $TX$; as for $\dn{\mu}{X}$, it
takes an argument of the form
\begin{displaymath}
\abst{s}{\pair{c[s]}{\abst{s'}{\pair{c'[s,s']}{f[s,s']}}}}  
\end{displaymath}
in $TTX$, and outputs
\begin{displaymath}
  \abst{s}{\pair{c'[s,c[s]]}{f[s,c[s]]}}
\end{displaymath}
which belongs to $TX$. In other words, $\dn{\mu}{X}$
is $\eval_{S\times X}^{S}$.
Recall that $T$-algebras are pairs $\pair{X}{h}$ where $X$
and $h:TX\to X$ are such that the following diagrams commute:
\begin{equation}
  \begin{xy}
    \xymatrix{T^2X\ar[r]^{Th}\ar[d]_{\mu_X} & TX\ar[d]^{h}\\
              TX\ar[r]_{h}  &  X}
  \end{xy}\qquad
  \begin{xy}
    \xymatrix{X\ar[r]^{\eta_X}\ar[dr]_{\id_X} & TX\ar[d]^{h} \\
                & X} 
  \end{xy}
\label{eq:alg}
\end{equation}
which means in this case that 
\begin{eqnarray}
h(\abst{s}{\pair{c[s]}{h(\abst{s'}{\pair{c'[s,s']}{x[s,s']}})}}) & = &
        h(\abst{s}{\pair{c'[s,c[s]]}{x[s,c[s]]}}) \label{eq:assoc}\\
h(\abst{s}{\pair{s}{x}}) & = & x \label{eq:ident}
\end{eqnarray}
A {\em morphism} of $T$-algebras $\pair{X}{h}$, $\pair{X'}{h'}$ is a map
$f:X\to X'$ making the following diagram commutative:
\begin{displaymath}
  \begin{xy}
    \xymatrix{ TX\ar[r]^{Tf}\ar[d]_h & TX'\ar[d]^{h'} \\
               X\ar[r]_f  & X'}
  \end{xy}
\end{displaymath}
Algebras and morphisms build a category $\ctg{Alg}^T$, with
an obvious forgetful functor $U^T:\pair{X}{h}\mapsto X$.

As with any adjunction, we may define a {\em comparison functor}
\begin{displaymath}
  K:\ctg{C}\to \ctg{Alg}^T
\end{displaymath}
which is given, in this case, by
\begin{displaymath}
  Y\mapsto \pair{Y^S}{\eval_Y^S}
\end{displaymath}

Call $U$ {\em monadic} if $K$ is an equivalence of categories.
Note that the same terminology is often applied to a stronger
notion, where $K$ is required to be an isomorphism
(see \cite{beck:tralco} and \cite{maclane:catfwm} on this issue).

The key idea is now that, in $\ctg{Sets}$, $T$-algebras are essentially
sets of maps $Y^S$, with their evaluation morphisms.
Theorem~\ref{thm:setmonad} below is a precise formulation of
this remark.
\begin{thm}\label{thm:setmonad}
  If $S\neq\emptyset$, then $U$ is monadic.
\end{thm}
\pre
Beck's criterion (\cite{beck:tralco}) applies, proving the statement.
However the situation is better understood by a direct
study of the comparison functor.

In order to show that $K$ is an equivalence,
we define a functor
\begin{displaymath}
  L:\ctg{Alg}^T\to \ctg{Sets}
\end{displaymath}
such that $KL$ and $LK$ are naturally isomorphic
to the identity functor  on $\ctg{Alg}^T$ and $\ctg{Sets}$
respectively.
Let $\pair{X}{h}$ be a $T$-algebra, we define the set
\begin{displaymath}
  Y=L\pair{X}{h}
\end{displaymath}
as follows.
From the projection
\begin{displaymath}
  \dn{q}{S\times X}: S\times (S\times X)\to S\times X
\end{displaymath}
we get
\begin{displaymath}
  \x{q}{S\times X}{*}: S\times X \to (S\times X)^S=TX
\end{displaymath}
hence
\begin{displaymath}
 \phi = h\cb \x{q}{S\times X}{*}: S\times X\to X
\end{displaymath}
Then $Y=L\pair{X}{h}$ will be the image of $\phi$. 
Concretely, $Y$ is the subset of
$X$ given by
\begin{displaymath}
  Y=\ens{h(\abst{s'}{\pair{s}{x}})|s\in S, x\in X}
\end{displaymath}
and we get a surjective map 
$e_X:S\times X\to Y$ making the following diagram commutative:
\begin{displaymath}
  \begin{xy}
    \xymatrix{S\times X\ar[d]_{e_X}\ar[r]^{\x{q}{S\times X}{*}} & TX\ar[d]^h \\
               Y\ar[r]_{m_X} & X }
  \end{xy}
\end{displaymath}
where $m_X$ is the inclusion monomorphism (in fact $e_X$ and $m_X$ 
really depend on the {\em algebra} $\pair{X}{h}$). 
Now $L$ has to be defined on morphisms as well: let then $u$ be a morphism
of algebras $u:\pair{X}{h}\to\pair{X'}{h'}$, $Y=L\pair{X}{h}$ and
$Y'=L\pair{X'}{h'}$. There is a unique
map $Lu:Y\to Y'$ such that the following diagram commutes:
\begin{displaymath}
  \begin{xy}
    \xymatrix{S\times X \ar[r]^{S\times u}\ar[d]_{e_X}&
   S\times X'\ar[d]^{e_{X'}}\\
              Y\ar[r]_{Lu} & Y'}
  \end{xy}
\end{displaymath}
Display $e_X$ as a cokernel
\begin{displaymath}
  \begin{xy}
    \xymatrix{Z\doubr{a}{b} & S\times X\ar[r]^{e_X} & Y}
  \end{xy}
\end{displaymath}
and fit the above diagram into this larger one:
\begin{displaymath}
  \begin{xy}
    \xymatrix{Z\doubr{a}{b} & S\times X \ar[ddd]_{e_X}\ar[rrr]^{S\times u}
            \ar[dr]^{\x{q}{S\times X}{*}}
             & & & S\times X'\ar[ddd]^{e_{X'}}\ar[dl]_{\x{q}{S\times X'}{*}}\\
                &           & TX\ar[r]^{Tu}\ar[d]_h  & TX'\ar[d]^{h'}   &  \\
                &           & X\ar[r]_u  & X'    &            \\
                & Y \ar[ur]_{m_X}\ar[rrr]_{Lu}&  &   & Y'\ar[ul]^{m_{X'}}             }
  \end{xy}
\end{displaymath}
The inner square commutes because $u$ is a morphism of algebras,
the top trapezoid commutes by naturality of $\up q*$, and the left
and right hand trapezoids commute
by definition of $e$, $m$ . Hence
\begin{eqnarray*}
  m_{X'}\cb e_{X'}\cb (S\times u)\cb a & = & u\cb m_X\cb e_X\cb a \\
                   & = & u\cb m_X\cb e_X\cb b \\
                   & = & m_{X'}\cb e_{X'}\cb (S\times u)\cb b
\end{eqnarray*}
and $m_{X'}$ being a monomorphism,
\begin{displaymath}
  e_{X'}\cb (S\times u)\cb a = e_{X'}\cb(S\times u)\cb b
\end{displaymath}
so that $e_{X'}\cb (S\times u)$ equalizes the pair $a$, $b$.
As $e_X$ is a cokernel,
there is a unique $Lu$ making the outer square commutative, as required.
As usual, uniqueness of $Lu$ ensures functoriality. 
Note that $e$ is a natural transformation:
\begin{displaymath}
  e:FU^T\to L
\end{displaymath}
As a consequence, the adjunction gives
a natural transformation from $U^T$ to $UL$, or better, as $U=U^TK$
\begin{displaymath}
 \up  e*: U^T\to U_TKL
\end{displaymath}
Thus we reduced the existence of a natural isomorphism:
\begin{displaymath}
  KL\simeq 1
\end{displaymath}
to the following statements:
\begin{itemize}
  \item for all $\pair{X}{h}$, $\x eX*$ is a bijection;

  \item for all $\pair{X}{h}$, $\x eX*$ is a morphism of algebras.
\end{itemize}
which will be  proved separately
in Lemmas~\ref{lem:iso} and \ref{lem:alg} below.

Finally, there is a natural isomorphism:
\begin{displaymath}
  LK\simeq 1
\end{displaymath}
If $Y$ is a set, $LKY$ is nothing but the subset of $Y^S$ consisting
of {\em constant maps} from $S$ to $Y$, which is naturally isomorphic
to $Y$ because $S\neq\emptyset$.
\erp

The two following lemmas will complete our argument. Keeping the 
previous notations,
\begin{lem}\label{lem:iso}
  For each $T$-algebra $\pair{X}{h}$, the map $\x eX*$ is a bijection.
\end{lem}
\pre
Let us give an explicit construction of the inverse map $f$.
For each $y\in\up{Y}{S}$, we define $f(y)\in X$
by
\begin{displaymath}
  f(y)=h(\abst{s}{\pair{s}{y[s]}})
\end{displaymath}
Now, for each $x\in X$,
\begin{eqnarray*}
  f(\x eX*(x)) & = & h(\abst{s}{\pair{s}{h(\abst{s'}{\pair{s}{x}})}})\\
               & = & h(\abst{s}{\pair{s}{x}}) \ \mathrm{by\ (\ref{eq:assoc})}\\
               & = & x  \ \mathrm{by\ (\ref{eq:ident})}
\end{eqnarray*}
Thus $f$ is a retraction for $\x eX*$.

On the other hand, let $y\in \up YS$,
\begin{eqnarray*}
  \x eX*(f(y)) & = & \abst{s}{h(\abst{s'}{\pair{s}{f(y)}})} \\
               & = & \abst{s}{h(\abst{s'}{\pair{s}
                     {h(\abst{s''}{\pair{s''}{y[s'']}})}}} \\
               & = & \abst{s}{h(\abst{s'}{\pair{s}{y[s]}})} \ 
                      \mathrm{by\ (\ref{eq:assoc})}
\end{eqnarray*}
but $y[s]=h(\abst{s''}{\pair{c[s]}{x[s]}})$ for suitable maps $c$ and $x$,
by definition of $Y$, so that the last expression reduces to
\begin{displaymath}
  \abst{s}{h(\abst{s'}{\pair{c[s]}{x[s]}})}=\abst{s}{y[s]}=y
\end{displaymath}
by (\ref{eq:assoc}) again.
Thus $\x eX*(f(y))=y$, and $f$ is also a section of $\x eX*$.
Hence $f=(\x eX*)^{-1}$, as required.
\erp

\begin{lem}\label{lem:alg}
  For each $T$-algebra $\pair{X}{h}$, $\x eX*$ is a morphism
  of algebras.
\end{lem}
\pre  Let us prove the commutativity of the following
diagram:
\begin{displaymath}
  \begin{xy}
    \xymatrix{TX\ar[r]^{T\x eX*}\ar[d]_h& T(Y^S)\ar[d]^{\eval_Y^S}\\
              X\ar[r]_{\x eX*}  & Y^S}
  \end{xy}
\end{displaymath}
Let $u\in TX$. It is a $\abst{s}{\pair{c[s]}{x[s]}}$ for a certain
pair of maps $c$, $x$. Hence 
\begin{eqnarray*}
  \x eX*\cb h(u) & = & \abst{s'}{h(\abst{s''}{\pair{s'}
                  {h(\abst{s}{\pair{c[s]}{x[s]}})}})}\\
             & = & \abst{s'}{h(\abst{s''}{\pair{c[s']}{x[s']}})}\ 
                  \mathrm{by\ (\ref{eq:assoc})}
\end{eqnarray*}
On the other hand
\begin{eqnarray*}
  \eval_Y^S\cb T\x eX*(u) & = & \eval_Y^S
      (\abst{s}{\pair{c[s]}{\abst{s'}{h(\abst{s''}{\pair{s'}{x[s]}})}}}) \\
                      & = & \abst{s}{h(\abst{s''}{\pair{c[s]}{x[s]}})}
\end{eqnarray*}
so that
\begin{displaymath}
  \x eX*\cb h = \eval_Y^S\cb T\x eX*
\end{displaymath}
\erp

\section{A general theorem}\label{sec:general}

By carefully examining the proofs in section \ref{sec:sets},
we see that the equations
(\ref{eq:assoc}) and (\ref{eq:ident}) are only used in
proving Lemmas~\ref{lem:iso} and~\ref{lem:alg}. Thus
large parts of our argument still hold beyond $\ctg{Sets}$.
As we shall see, a key hypothesis for generalizing
the previous results is that
$\ctg{C}$ {\em has regular epi-mono factorization},
that is, each arrow in $\ctg{C}$ factorizes
as a composition of a regular epimorphism and a monomorphism.
The regularity hypothesis implies that such a factorization is
unique, up to isomorphism. 
Thus from now on, we assume that $\ctg{C}$ has regular epi-mono factorization.

\subsection{Construction of $L$}\label{subsec:constrL}

Let us first generalize the
construction of
\begin{displaymath}
  L:\ctg{Alg}^T\to \ctg{C}
\end{displaymath}
If $\ctg{C}$ is any cartesian closed category, and
$\pair{X}{h}$ is a $T$-algebra, the morphism $f_X=h\cb\x q{S\times X}*$
is still defined:
\begin{displaymath}
  \begin{xy}
    \xymatrix{S\times X\ar[r]^{\x q{S\times X}*}& TX\ar[r]^h & X}
  \end{xy}
\end{displaymath}
so that all we need to extend the previous definition
of $L$ is the existence of an image object for $f_X$, which immediately
follows from the regular epi-mono factorization property on $\ctg{C}$.

Thus we get
a regular epimorphism $e_X$ and a monomorphism $m_X$ such that
\begin{equation}
  f_X=m_X\cb e_X
\label{eq:epimono}
\end{equation}
and $L\pair{X}{h}$ can be defined as the object $Y$ in the commutative
diagram
\begin{displaymath}
  \begin{xy}
    \xymatrix{S\times X\ar[r]^{q_{S\times X}^*}\ar[d]_{e_X} & TX \ar[d]^h\\
                 Y\ar[r]_{m_X}   & X}
  \end{xy}
\end{displaymath}
Now the regularity of $e_X$ suffices to make $L$ a functor,
and $e$ a natural transformation:
\begin{displaymath}
  e:FU^T\to L
\end{displaymath}
as already shown in the proof of Theorem~\ref{thm:setmonad}.
Whence a natural transformation
\begin{displaymath}
  \up e*:U^T\to U^TKL
\end{displaymath}

\subsection{A retraction}

\begin{prop}\label{prop:retrac}
  For each $T$-algebra $\pair{X}{h}$, $\x eX*:X\to Y^S$ has a retraction.
\end{prop}
\pre
Consider the following commutative diagram:
\begin{equation}
  \begin{xy}
    \xymatrix{S\times X\ar[r]^{\x q{S\times X}*}\ar[d]_{\dn eX}\ar[dr]^{\dn fX}
                              & TX\ar[d]^h \\
              Y\ar[r]_{\dn mX}         & X}
  \end{xy}
\label{eq:retrac}
\end{equation}
The left lower triangle gives rise by adjunction to the following
commutative triangle:
\begin{displaymath}
  \begin{xy}
    \xymatrix{X\ar[d]_{\x eX*}\ar[dr]^{\x fX*} & \\
              \up YS\ar[r]_{\x mXS} & \up XS}
  \end{xy}
\end{displaymath}
If $\x fX*$ has a retraction $r$, then $r\cb\x mXS$ is immediately
a retraction for $\x eX*$. Thus we turn to the construction 
of such a retraction for $\x fX*$. The upper right triangle
of (\ref{eq:retrac}) gives rise again by adjunction to 
the following commutative triangle:
\begin{equation}
  \begin{xy}
    \xymatrix{X\ar[r]^{(\x q{S\times X}*)^*}\ar[dr]_{\x fX*} & 
             \up{(TX)}S\ar[d]^{\up hS} \\
                & \up XS}
  \end{xy}
\label{eq:triangle}
\end{equation}
Let us introduce a natural transformation $\theta$ between $U$ and $T$.
For each object $Z$, we first have
\begin{displaymath}
  \pair{\dn p{Z^S}}{\eval_Z}:S\times \up ZS\to S\times Z
\end{displaymath}
whence
\begin{displaymath}
  \theta_Z=\pair{\dn p{Z^S}}{\eval_Z}^*:\up ZS\to TZ
\end{displaymath}
as required. The triangle (\ref{eq:triangle}) then fits into 
a larger diagram:
\begin{equation}
  \begin{xy}
    \xymatrix{
X\ar[r]^{(\x q{S\times X}*)^*}\ar[dr]_{\x fX*} & 
        (TX)^S\ar[r]^{\dn{\theta}{TX}}\ar[d]_{h^S} & 
        TTX\ar[d]^{Th} \\
  & X^S\ar[r]_{\dn{\theta}X}\ar[dr]_r & TX\ar[d]^h   \\
  &        & X }
  \end{xy}
\label{eq:bigdiagram}
\end{equation}
where $r$ is defined by
\begin{displaymath}
  r=h\cb\dn{\theta}X
\end{displaymath}
Now (\ref{eq:bigdiagram}) commutes: in fact (\ref{eq:triangle})
commutes already, the lower right triangle commutes by definition
of $r$, and the square by naturality of $\theta$. By defining
\begin{displaymath}
  \psi=\dn{\theta}{TX}\cb (\x q{S\times X}*)^*
\end{displaymath}
we get
\begin{displaymath}
  r\cb \x fX* = h\cb Th\cb\psi
\end{displaymath}
From the definition of a $T$-algebra~(\ref{eq:alg}), we know
that $h\cb Th=h\cb \mu_X$. We claim that
\begin{equation}
  \mu_X\cb\psi=\eta_X
\label{eq:muetapsi}
\end{equation}
which implies
\begin{eqnarray*}
  r\cb\x fX* & = & h\cb Th\cb\psi\\
          & = & h\cb\mu_X\cb\psi\\
          & = & h\cb\eta_X\\
          & = & \id_X
\end{eqnarray*}
by~(\ref{eq:alg}) again, and the result reduces to
(\ref{eq:muetapsi}), proved in Lemma~\ref{lem:muetapsi}.
\erp

\begin{lem}\label{lem:muetapsi}
  For each object $X$, $\mu_X\cb\theta_{TX}\cb(\x q{S\times X}*)^*=\eta_X$.
\end{lem}
\pre
The lemma states the commutativity of the following diagram:
\begin{equation}
  \begin{xy}
    \xymatrix{(TX)^S\ar[r]^{\dn\theta{TX}} & TTX\ar[d]^{\dn\mu{X}} \\
   X\ar[u]^{\up{(\x q{S\times X}*)}*}\ar[r]_{\dn\eta{X}}\ar[ur]^{\psi}& TX}
  \end{xy}
\label{eq:thetaetapsi}
\end{equation}
where the upper left triangle commutes by definition of $\psi$. The first
step is to note that, as
$\dn\theta{TX}=\up{\pair{\dn{p}{TX^S}}{\eval_{TX}}}*$, $\psi$ is
of the form $\up\phi*$, where $\phi$ makes the following triangle
commute:
\begin{equation}
  \begin{xy}
    \xymatrix{S\times (TX)^S
\ar[r]^*+<10pt>{{}^{\pair{\dn{p}{TX^S}}{\eval_{TX}}}}
         & S\times TX \\
       S\times X\ar[u]^{S\times\up{(\x q{S\times X}*)}*}\ar[ru]_{\phi} & }
  \end{xy}
\end{equation}
Let us show that
\begin{equation}
  \phi = \pair{\dn pX}{\x q{S\times X}*}
\label{eq:phi}
\end{equation}
This amounts to check the commutativity in the product diagram:
\begin{equation}
  \begin{xy}
    \xymatrix{ & S\times X\ar[ld]_{\dn pX}\ar[d]^{\phi}
      \ar[rd]^{\x q{S\times X}*} & \\
              S & S\times TX\ar[l]^{\dn p{TX}}\ar[r]_{\dn q{TX}} & TX}
  \end{xy}
\end{equation}
The left hand side commutes because of the commutativity of:
\begin{equation}
  \begin{xy}
    \xymatrix{S\times X\ar[r]^*+<10pt>{{}^{S\times\up{(\x q{S\times X}*)}*}}
\ar[dr]_{\dn pX} & 
S\times (TX)^S\ar[r]^*+<10pt>{{}^{\pair{\dn p{TX^S}}{\eval_{TX}}}}\ar[d]|{\dn p{TX^S}} & 
   S\times TX\ar[dl]^{\dn p{TX}} \\
                        & S              &         } 
  \end{xy}
\end{equation}
As for the right hand side, it amounts to the commutativity of:
\begin{equation}
  \begin{xy}
    \xymatrix{S\times X\ar[r]^*+<10pt>{{}^{S\times\up{(\x q{S\times X}*)}*}}
\ar[dr]_{\x q{S\times X}*} & 
S\times (TX)^S\ar[r]^*+<10pt>{{}^{\pair{\dn p{TX^S}}{\eval_{TX}}}}\ar[d]|{\eval_{TX}} & 
   S\times TX\ar[dl]^{\dn q{TX}} \\
                        & TX             &         } 
  \end{xy}
\label{eq:qcommut}
\end{equation}
Now, in~(\ref{eq:qcommut}), the commutativity of the right hand side
is straightforward, so we are reduced to check commutativity on
the left hand side, which in turn results from:
\begin{eqnarray*}
  \eval_{TX}\cb (S\times\up{(\x q{S\times X}*)}*) & = 
             & {(\id_{\up{(TX)}S})}_*\cb (S\times\up{(\x q{S\times X}*)}*) \\
   & = & ({\id_{\up{(TX)}S}}\cb \up{(\x q{S\times X}*)}*)_*\\
   & = & {(\up{(\x q{S\times X}*)}*)}_*\\
   & = & \x q{S\times X}*
\end{eqnarray*}
This achieves the proof of~(\ref{eq:phi}).

Back to (\ref{eq:thetaetapsi}),
we now prove that the right lower triangle commutes; we just proved
that this triangle is in fact:
\begin{equation}
  \begin{xy}
    \xymatrix{ & TTX\ar[d]^{\mu_X}\\
              X\ar[ur]^{\pair{\dn pX}{\x q{S\times X}*}^*}\ar[r]_{\eta_X}& TX}
  \end{xy}
\label{eq:etamu}
\end{equation}
As $\mu_X=(\eval_{S\times X})^S$
and $\eta_X=(\id_{S\times X})^*$, (\ref{eq:etamu}) commutes
if and only if~(\ref{eq:ideval}) also commutes:
\begin{equation}
  \begin{xy}
    \xymatrix{ & S\times TX\ar[d]^{\eval_{S\times X}}\\
 S\times X\ar[ur]^{\pair{\dn pX}{\x q{S\times X}*}}\ar[r]_{\id_{S\times X}}
 & S\times X}
  \end{xy}
\label{eq:ideval}
\end{equation}
Let $\Delta=\pair{\id_S}{\id_S}:S\to S\times S$, we may express 
$\pair{\dn pX}{\x q{S\times X}*}$ as a compostion:
\begin{equation}
  \pair{\dn pX}{\x q{S\times X}*}=(S\times\x q{S\times X}*)\cb(\Delta\times X)
\label{eq:compos}
\end{equation}
such that the commutativity of~(\ref{eq:ideval}) reduces to the commutativity
of the following diagram:
\begin{equation}
  \begin{xy}
    \xymatrix{S\times S\times X\ar[r]^*+<10pt>{{}^{S\times\x q{S\times X}*}}
        & S\times TX\ar[d]^{\eval_{S\times X}}\\
                S\times X\ar[r]_{\id_{S\times X}}\ar[u]^{S\times\Delta}
         & S\times X}
  \end{xy}
\end{equation}
but the adjunction gives:
\begin{eqnarray*}
  \eval_{S\times X}\cb (S\times\x q{S\times X}*) & = 
                & (\id_{TX})_*\cb (S\times\x q{S\times X}*)\\
               & = & ({\id_{TX}}\cb{\x q{S\times X}*})_*\\
               & = & (\x q{S\times X}*)_* \\
               & = & q_{S\times X}
\end{eqnarray*}
Finally
\begin{displaymath}
  q_{S\times X}\cb(\Delta\times X) = \id_{S\times X}
\end{displaymath}
and the lemma is proved.
\erp

\subsection{Existence of a section}
We now turn to conditions ensuring the existence of a section
for $\x eX*$. We first prove a technical lemma:
\begin{lem}\label{lem:esect}
  Suppose $\pair{X}{h}$ is a $T$-algebra, and the epimorphism
  $\dn eX$ has a section. Then $\x eX*$ also has a section.
\end{lem}
\pre
Let $\sigma$ be a section of $\dn eX$, that is
\begin{equation}
  \dn eX\cb\sigma=\id_Y
\end{equation}
and define
\begin{equation}
  \Sigma=h\cb\up\sigma S
\label{eq:section}
\end{equation}
we claim that $\Sigma$ is a section of $\x eX*$. Consider
the following diagram:
\begin{equation}
  \begin{xy}
    \xymatrix{Y^S\ar[r]^{\up\sigma{S}}\ar[dr]_{\Sigma} &
          TX\ar[r]^{\up{(\x q{S\times TX}*)}*}\ar[d]_h &
           (TTX)^S\doubd{\up{(Th)}S}{\x{\mu}XS} \\
                  & X\ar[r]_{\up{(\x q{S\times X}*)}*}\ar[d]_{\x eX*} & 
            (TX)^S\ar[d]^{\up hS}   \\
                  & Y^S\ar[r]_{\x mXS}  & X^S   }
  \end{xy}
\end{equation}
and note that
\begin{eqnarray*}
  \x mXS\cb\x eX*\cb\Sigma & = &
          \x mXS\cb\x eX*\cb h\cb\up{\sigma}S \\
    & = & \up hS\cb\up{(\x q{S\times X}*)}*\cb h\cb\up{\sigma}S \\
    & = & \up hS\cb\up{(Th)}S\cb\up{(\x q{S\times TX}*)}*\cb\up{\sigma}S\\
    & = & \up hS\cb\x{\mu}XS\cb\up{(\x q{S\times TX}*)}*\cb\up{\sigma}S 
\end{eqnarray*}
but using the adjunctions:
\begin{eqnarray*}
 \x{\mu}XS\cb\up{(\x q{S\times TX}*)}{*}& = 
         & {(\dn{\mu}X\cb \x q{S\times TX}*)}^{*} \\
          & = & {(\x\eval{S\times X}S\cb\x q{S\times TX}*)}^*\\
          & = & {(\x q{S\times X}*\cb\dn{\eval}{S\times X})}^*\\
          & = & {(\x q{S\times X}*)}^S\cb {\x{\eval}{S\times X}*}\\
          & = & {(\x q{S\times X}*)}^S\cb {((\id_{TX})_*)}^*\\
          & = & {(\x q{S\times X}*)}^S
\end{eqnarray*}
so that
\begin{eqnarray*}
  \x mXS\cb \x eX*\cb\Sigma & = & \up hS\cb {(\x q{S\times X}*)}^S
                                       \cb\up{\sigma}S\\
                     & = & \x mXS\cb\x eXS\cb\up{\sigma}S\\
                     & = & \x mXS
\end{eqnarray*}
and finally
\begin{equation}
  \x eX*\cb\Sigma = \id_{Y^S}
\end{equation}
because $\x mXS$ is a monomorphism (right-adjoints preserve monomorphisms).
\erp

Now $e_X$ was defined as a regular epimorphism,
which does not imply the existence of a section:
think for example at the presheaf category of graphs, where all
epimorphisms are regular, many of them without a section.
Nevertheless, a very simple condition on $S$ will be sufficient:
recall that a (global) {\em element} of an object $Z$ in $\ctg{C}$
is an arrow $z:1\to Z$, where $1$ is the terminal object,
then
\begin{lem}\label{lem:existsec}
  Suppose $S$ has at least one element and $\pair{X}{h}$
  is a $T$-algebra. Then $\dn eX$ has a section. 
\end{lem}
\pre
Let $s_0:1\to S$ be an element of $S$: for each object $Z$, it
gives rise to 
\begin{displaymath}
  \gamma_Z:Z^S\to Z
\end{displaymath}
by composition of $Z^{s_0}$ with the canonical isomorphism $A^1\simeq A$.
Note that $\gamma$ becomes a natural transformation and that for each $Z$
\begin{equation}
  \gamma_Z\cb \x qZ* =\id_Z
\label{eq:gq}
\end{equation}
Let us now consider an algebra $\pair{X}{h}$, and build the following
diagram:
\begin{displaymath}
  \begin{xy}
    \xymatrix{S\times TX\ar@<-2pt>[r]_{\x q{S\times TX}*}
     \ar[d]_{S\times h}   & 
      TTX\ar@<-2pt>[d]_{Th}\ar@<-2pt>[l]_{\gamma_{S\times TX}} \\
S\times X\ar@<-2pt>[r]_{\x q{S\times X}*}
          \ar[d]_{\dn eX}    & TX\ar@<-2pt>[l]_{\gamma_{S\times X}}
            \ar@<-2pt>[d]_h\ar@<-2pt>[u]_{\eta_{TX}}  \\
Y\ar[r]_{\dn mX} & X\ar@<-2pt>[u]_{\eta_X}}
  \end{xy}
\end{displaymath}
Define
\begin{equation}
  \sigma=\gamma_{S\times X}\cb\eta_X\cb m_X
\end{equation}
we claim that $\sigma$ is a section of $e$. To see this, we compute
\begin{equation}
  m_X\cb e_X\cb \sigma\cb e_X
\label{eq:mesigma}
\end{equation}
First, using naturality of $\eta$, $\gamma$ and $\up q*$:
\begin{eqnarray*}
  m_X\cb e_X\cb \sigma\cb e_X & = 
 & h\cb \x q{S\times X}*\cb\dn \gamma{S\times X}\cb\dn{\eta}X\cb h
        \cb\x q{S\times X}*\\
 & = & 
h\cb\x q{S\times X}*\cb\dn \gamma{S\times X}\cb Th\cb\dn{\eta}{TX}
\cb\x q{S\times X}*\\
&=&
h\cb\x q{S\times X}*\cb (S\times h)\cb\dn{\gamma}{S\times TX}
\cb\dn{\eta}{TX}\cb\x q{S\times X}*\\
&=&
h\cb Th\cb\x q{S\times TX}*\cb\dn{\gamma}{S\times TX}
\cb\dn{\eta}{TX}\cb\x q{S\times X}*
\end{eqnarray*}
then because $h\cb Th=h\cb\mu_X$ and $\mu_X=\x{\eval}{TX}S$,
\begin{displaymath}
 m_X\cb e_X\cb \sigma\cb e_X =
h\cb \x{\eval}{TX}S\cb\x q{S\times TX}*\cb\dn{\gamma}{S\times TX}
\cb\dn{\eta}{TX}\cb\x q{S\times X}*
\end{displaymath}
and by using naturality of $\dn q*$ and $\gamma$ again, (\ref{eq:mesigma})
reduces to:
\begin{equation}
  h\cb \x q{S\times X}* \cb \dn{\gamma}{S\times X}\cb \x{\eval}{TX}S
    \cb \dn{\eta}{TX}\cb \x q{S\times X}*
\label{eq:reduc}
\end{equation}
Now
\begin{displaymath}
  \x{\eval}{TX}S\cb \dn{\eta}X = \dn{\mu}X\cb\dn{\eta}{TX} = \id_{TX}
\end{displaymath}
and 
\begin{displaymath}
  \dn{\gamma}{S\times X}\cb\x q{S\times X}* = \id_{S\times X}
\end{displaymath}
so that finally
\begin{eqnarray*}
 m_X\cb e_X\cb \sigma\cb e_X & = & h\cb \x q{S\times X}*\\
                             & = & m_X\cb e_X 
\end{eqnarray*}
As $m_X$ is a monomorphism, and $e_X$ an epimorphism, this implies
\begin{displaymath}
  e_X\cb \sigma = \id_Y
\end{displaymath}
which ends the proof.
\erp
From Lemmas~\ref{lem:esect} and~\ref{lem:existsec} immediately follows
\begin{prop}\label{prop:section}
  If $S$ has at least one element, then for each $T$-algebra
$\pair{X}{h}$, $\x eX*$ has a section.
\end{prop}

\subsection{$\x eX*$ is a morphism of algebras}

We show that Lemma~\ref{lem:alg} extends to
the general setting without additional assumptions:

\begin{lem}\label{lem:algbis}
  For each $T$-algebra $\pair{X}{h}$, $\x eX*$ is a morphism
of algebras.
\end{lem}
\pre
Let us consider the following diagram:
\begin{displaymath}
  \begin{xy}
    \xymatrix{TX\ar[r]^{T(\x eX*)}\ar[d]_h &
 T(\up YS)\ar[r]^{T(\x mXS)}\ar[d]^{\x{\eval}YS} & 
 T(\up XS)\ar[d]^{\x{\eval}XS}\\
              X\ar[r]_{\x eX*} & \up YS\ar[r]_{\x mXS} & \up XS}
  \end{xy}
\end{displaymath}
We must show that the left hand side square commutes.
First, by using the adjunction, together with (\ref{eq:alg}),
(\ref{eq:epimono}) and the fact
that $\up q*$ is natural,
\begin{eqnarray*}
  \x mXS\cb \x eX* \cb h & = & \up{(\dn mX\cb\dn eX)}*\cb h\\
                         & = & \up{(\dn mX\cb \dn eX\cb (S\times h))}*\\
                         & = & \up{(h\cb \x q{S\times X}*\cb (S\times h))}*\\
                         & = & \up{(h\cb Th \cb \x q{S\times TX}*)}*\\
                         & = & \up{(h\cb \mu_X\cb \x q{S\times TX}*)}*\\
                         & = & \up{(h\cb \x{\eval}{S\times X}S
                                                \cb \x q{S\times TX}*)}*\\
                         & = & \up{(h\cb \x q{S\times X}*\cb 
                                          \dn{\eval}{S\times X})}*
\end{eqnarray*}
then the same ingredients plus the facts that $T$ is a functor, and 
$\eval^S$ is natural imply:
\begin{eqnarray*}
  \x mXS\cb \x{\eval}YS\cb T(\x eX*) & = & \x{\eval}XS\cb T(\x mXS)
                                             \cb T(\x eX*)\\
                      & = & \x{\eval}XS\cb T(\x mXS\cb \x eX*)\\
                      & = & \x{\eval}XS\cb T(\up{(h\cb\x q{S\times X}*)}*)\\
                      & = & \x{\eval}XS
                           \cb \up{(S\times\up{(h\cb\x q{S\times X}*)}*)}S\\
            & = & \up{(\eval_X\cb (S\times\up{(h\cb\x q{S\times X}*)}*))}S\\
            & = & \up{\dn{(\up{(h\cb \x q{S\times X}*)}*)}*}S\\
            & = & \up{(h\cb\x q{S\times X}*)}S\\
            & = & \up{(h\cb\x q{S\times X}*)}S\cb\up{(\dn{(\id_{TX})}*)}*\\
            & = & \up{(h \cb \x q{S\times X}*\cb \dn{(\id_{TX})}*)}*\\
            & = & \up{(h\cb \x q{S\times X}*\cb \dn{\eval}{S\times X})}*
\end{eqnarray*}
whence
\begin{displaymath}
  \x mXS\cb \x eX* \cb h = \x mXS\cb \x{\eval}YS\cb T(\x eX*)
\end{displaymath}
and because $\x mXS$ is a monomorphism,
\begin{displaymath}
  \x eX*\cb h = \x{\eval}YS\cb T(\x eX*)
\end{displaymath}
which ends the proof.
\erp

\subsection{Monadicity theorem}

As a consequence of propositions~\ref{prop:retrac} and~\ref{prop:section},
if $S$ has at least one element,
then $KL$ is naturally isomorphic to the identity on $\ctg{Alg}_T$.
On the other hand,
\begin{prop}\label{prop:lk}
  If $S$ has at least one element, $LK\simeq 1$.
\end{prop}
\pre
Let $Y$ be an object of $\ctg{C}$. The algebra $KY$ is
$\pair{Y^S}{\x{\eval}YS}$. By naturality of $\up q*$, the
following diagram commutes:
\begin{displaymath}
  \begin{xy}
    \xymatrix{S\times\up YS\ar[r]^{\x q{S\times{\up YS}}*}\ar[d]_{\eval_Y}
                & T(\up YS)\ar[d]^{\x{\eval}YS}\\
               Y\ar[r]_{\x qY*} & \up YS
}
  \end{xy}
\end{displaymath}
Now, $S$ has at least one element, say
\begin{displaymath}
  s_0:1\to S
\end{displaymath}
which easily provides a retraction for $\x qY*$, and a section for
$\dn{\eval}Y$. In particular, $\x qY*$ is a monomorphism and $\dn{\eval}Y$
a split --- thus regular --- epimorphism. Hence  there is a unique
isomorphism $\xi_Y$ making the following diagram commutative:
\begin{displaymath}
  \begin{xy}
    \xymatrix{S\times\up YS\ar[r]^{\dn e{Y^S}}\ar[d]_{\eval_Y}
                & LKY\ar[d]^{\dn m{Y^S}}\\
               Y\ar[r]_{\x qY*}\ar[ur]|{\xi_Y} & \up YS
}
  \end{xy}
\end{displaymath}
Moreover, $\xi_Y$ is natural in $Y$, as shown by standard uniqueness 
arguments.
\erp
We finally state our main result, an immediate consequence
of the above discussion:
\begin{thm}\label{thm:monadicity}
  If $\ctg{C}$ has regular epi-mono factorization,
  and $S$ is an object of $\ctg{C}$ having
  at least one element, then $U:X\mapsto X^S$ is monadic.
\end{thm}

\section{Equations}

As we pointed out in the introduction, the present
work is strongly related to the analysis of
global states given in \cite{plotkinpower:notcdm},
with a significant difference: in \cite{plotkinpower:notcdm},
$\ctg{C}$ is any category with countable products and coproducts,
and $S$ is a countable {\em set} of states, so that 
$S\times X$ (resp.\ $X^S$) now denotes the coproduct (resp.\ the product) 
of $S$ copies of $X$, and not as in cartesian closed categories
the internal product (resp.\ exponential) by an object $S$. 
Suppose however that $\ctg{C}$
satisfies both sets of conditions, those in \cite{plotkinpower:notcdm}
and those of the present paper, as for example will be the case of
any presheaf category:
then a countable set $S$ can be embedded
as an object of $\ctg{C}$, coproduct of $S$ copies of the
terminal element, and the results of \cite{plotkinpower:notcdm}
still make sense in our setting, and the two possible interpretations
of the notations $S\times -$ and $(-)^S$ coincide. 

In particular, under these additional hypotheses,
we may revisit the equational presentation
of global state algebras. Let $\Sigma$
be the signature consisting of
a symbol ${\term l}$ of arity $S$, and for each
$s\in S$, a symbol ${\term u}_s$, of arity $1$; let 
$E$ be the following set of equations among terms
generated by $\Sigma$:
\begin{eqnarray*}
  {\term u}_s[{\term u}_t[x]] & = & {\term u}_t[x]\\
  {\term u}_s[{\term l}[(a_t)_t]] & = & {\term u}_s[a_s]\\
  {\term l}[({\term u}_s[x])_s] & = & x\\
  {\term l}[({\term l}[(a_{st})_t])_s] & = & {\term l}[(a_{ss})_s]
\end{eqnarray*}
A $\pair{\Sigma}{E}$-algebra in $\ctg{C}$ is now a pair 
$\pair{A}{\abs{.}_A}$ where $A$ is an object of $\ctg{C}$,
and $\abs{.}_A$ assigns to each symbol in $\Sigma$ an arrow
of $\ctg{C}$ of appropriate arity:
\begin{eqnarray*}
  \abs{{\term l}}_A & : & A^S\to A\\
  \abs{{\term u}_s}_A & : & A\to A 
\end{eqnarray*}
in such a way that the equations in $E$ are satisfied.

Examples of $\pair{\Sigma}{E}$-algebras are objects of
the form $B^S$, with the following interpretation of 
$\Sigma$: $\abs{{\term l}}_A$ is defined by
\begin{displaymath}
  \begin{xy}
    \xymatrix{(B^S)^S\ar[r]^{\simeq} & B^{S\times S}\ar[r]^{B^\delta} & B^S}
  \end{xy}
\end{displaymath}
where $\delta:S\to S\times S$ is the diagonal map, and 
$\abs{{\term u}_s}_A$ by the composite:
\begin{displaymath}
  \begin{xy}
    \xymatrix{B^S\ar[r]^{B^s} & B\ar[r]^{B^!} & B^S}
  \end{xy}
\end{displaymath}
where $s:1\to S$ picks $s\in S$ and $!:S\to 1$. In $\ctg{Sets}$,
${\term l}$ takes a family of maps $(b_s)_{s\in S}$ and returns
the diagonal map $s\mapsto b_s(s)$, whereas ${\term u}_s$
takes a map $b$ and returns the constant map $s'\mapsto b(s)$.
Now, joigning Theorem~\ref{thm:monadicity} above and Theorem~1 of
\cite{plotkinpower:notcdm}, we may conclude that those are
essentially the {\em only} examples of $\pair{\Sigma}{E}$-algebras,
and of course of $T$-algebras.

\subsection*{Remark} Because we deal with global states only,
we take here the set of locations $L$ as a singleton,
which reduces the seven equations of \cite{plotkinpower:notcdm}, section~3
to the first four.

\subsubsection*{Acknowledgement} Many thanks to Albert Burroni
for numerous illuminating conversations on monadicity.

\end{document}